\documentclass[11pt,reqno]{amsart}
\usepackage{amsmath,amsfonts,amssymb,amsthm,amscd}
 
\makeindex

\usepackage{amssymb}
\usepackage[latin1]{inputenc}
\usepackage{ epsfig,epsf}
\usepackage{enumerate}
\usepackage{indentfirst}
\usepackage{euscript}
\usepackage{graphicx}
\usepackage{amsmath}
\usepackage{amsfonts}
\usepackage{amssymb}
\makeindex

 \theoremstyle{plain}
 
\newtheorem{proposition}{Proposition}
\newtheorem{lemma}{Lemma}

\newtheorem{remark}{Remark}
\theoremstyle{definition}

\newtheorem{theorem}{Theorem}

 \newcommand{\ep}{\epsilon}

\usepackage{amsmath,amsfonts,  amssymb,amsthm,amscd}
\makeindex
\usepackage{amssymb}
\usepackage{graphics}
\usepackage[latin1]{inputenc}
\usepackage{epsf}

\usepackage{ epsfig,epsf}
\usepackage{enumerate}
\usepackage{indentfirst}
\usepackage{euscript}
\usepackage{graphicx}
\usepackage{amsmath}
\usepackage{amsfonts}
\usepackage{amssymb}
\makeindex

 \theoremstyle{plain}

 \title[   Harmonic Mean Curvature Lines ]{   Harmonic Mean Curvature
Lines  on Surfaces Immersed  in $\mathbb R^3$ }

 \author[R. Garcia]{Ronaldo Garcia}

\author[J. Sotomayor]{Jorge Sotomayor}


 \thanks{The first author was partially supported by     FUNAPE/UFG.
Both   authors are fellows of  CNPq.
 This work was done under the project PRONEX/FINEP/MCT - Conv.
76.97.1080.00 - Teoria Qualitativa das Equa\c c\~oes Diferenciais
Ordin\'arias and CNPq - Grant 476886/2001-5.}

\begin{document}
 \maketitle

 \small {{\sc Abstract.-}
 Consider  oriented surfaces immersed in
$\mathbb R^3. $ Associated to them, here are studied  pairs of
transversal  foliations
 with singularities, defined on  the  {\it Elliptic} region, where the
Gaussian curvature  $\mathcal K$, given by the product of the
principal curvatures ${k_1}, k_2$ is positive.
 The   leaves of the foliations are the  {\it lines of    harmonic mean
curvature,} also called {\it characteristic } or {\it diagonal lines},
along which   the normal  curvature of the immersion is given by
${\mathcal K}/{\mathcal H}$,  where $ {\mathcal H}=({k_1}+k_2)/2$ is
the
arithmetic mean curvature.  That is, ${\mathcal K}/{\mathcal
H}=((1/{k_1} + 1/{k_2})/2)^{-1}$ is the  {\it    harmonic mean} of the
principal curvatures ${k_1}, k_2$ of the immersion.
The singularities of the foliations  are
 the {\it umbilic points} and {\it parabolic curves}, where ${k_1} =
k_2$ and  ${\mathcal K} = 0$, respectively.

Here are determined the structurally  stable patterns of     {\it
harmonic mean curvature lines}  near the
 {\it umbilic points},   {\it parabolic curves}  and {\it     harmonic
mean  curvature  cycles},   the periodic  leaves of the foliations. The
genericity of these patterns is established.

This  provides  the three essential local ingredients to establish
 sufficient
 conditions, likely to be also  necessary,  for {\it   Harmonic  Mean
Curvature Structural Stability} of immersed surfaces. This study,
outlined towards the end of the paper, is   a natural   analog
and
complement  for that carried out previously by the authors for the
{\it  Arithmetic Mean Curvature}  and  the {\it Asymptotic Structural
Stability} of immersed surfaces,  \cite {a1, m, a2}, and also extended
recently to the case of the {\it  Geometric Mean Curvature Configuration}
\cite {g}.}
\vskip .4cm

\section{Introduction}

In this  paper are  studied  the {\it   harmonic mean curvature
configurations} associated to  immersions of  oriented surfaces
into $\mathbb R^3$. They consist on the {\it umbilic points} and
{\it parabolic curves}, as singularities,   and of  the {\it lines
of  harmonic mean  curvature}  of the immersions, as the leaves of
the two transversal foliations in the configurations. The  normal
curvature of the immersion along these
lines    is  given by the {\it harmonic mean}  of
the principal curvatures, defined  by
${\mathcal K}/{\mathcal H}=((1/{k_1} + 1/{k_2})/2)^{-1}$, in
terms of the standard curvature functions:  {\it principal curvatures}
${k_1}, {k_2}$,  {\it   arithmetic mean curvature}  ${\mathcal H}=({k_1}
+ {k_2})/2$ and {\it   Gaussian  curvature}  ${\mathcal K}=
{k_1}{k_2}$.

The two transversal foliations, called here {\it   harmonic mean
curvature foliations},
are well defined and regular only on the non-umbilic part of the
elliptic region of the
immersion, where the  Gaussian Curvature is positive. In fact, there
they are the integral curves
of smooth quadratic differential equations.  The set where the Gaussian
Curvature vanishes,
 the parabolic set, is generically a  regular curve which is the border
of  the elliptic region; see \cite {bgm}.
The umbilic points
 are those at which the principal curvatures coincide,  generically are
isolated and disjoint
 from the parabolic curve. See section 2 for precise definitions.

This study is  a natural development and extension of  previous results
about  the
 Arithmetic Mean Curvature and Asymptotic
Configurations,
 dealing with the qualitative properties of the  lines along which the
normal  curvature
is  the  arithmetic mean of the principal curvatures (i.e. is  the
standard Mean Curvature) or is null.
This  has been considered  previously by the authors; see  \cite {a1,
a2} and  \cite {m}, and has  also been extended recently to the case of
the {\it  Geometric Mean Curvature} \cite {g}.

The point of departure of this line of research, however, can be found
in the classical works of Euler,  Monge, Dupin and Darboux, concerned
with the  lines of principal curvature and umbilic points  of
immersions. See \cite {ca, Sp, St} for an initiation on  the basic facts on this
subject; see  \cite {gs1, gs2} for a  discussion of the classical
contributions and for their analysis  from the point of view of  structural
stability  of differential equations. A modern general presentation of
structural stability of dynamical systems can be found in  \cite {pm}.

This paper establishes  sufficient conditions, likely to be also
necessary, for the  structural stability of {\it   harmonic mean curvature
configurations}
  under small perturbations of the immersion. See section 7  for
precise statements.

  This  extends to the   harmonic mean curvature setting the  main
theorems on structural stability   for the arithmetic and  geometric mean
curvature
configurations and for the asymptotic configurations, proved in \cite
{a1, m, g, a2}.

Three  local ingredients  are essential  for this extension: the
umbilic
points, endowed  with  their  harmonic  mean curvature separatrix
structure, the  harmonic
mean curvature cycles, with the calculation of the derivative of the
Poincar\'e return map,  through which is expressed the hyperbolicity
condition and   the parabolic curve, together with the parabolic tangential
singularities
and associated separatrix structure.

The conclusions  of this paper, on the elliptic  region, are
complementary
to  results valid independently on the
hyperbolic region (on which the Gaussian curvature is negative), where
the separatrix  structure near the parabolic curve and the asymptotic
structural stability has been studied in \cite {a1, a2}.

The parallel  with  the conditions for principal, arithmetic  mean
curvature and asymptotic  structural stability is remarkable. This can be
attributed to the unifying role played by the notion of Structural
Stability  of  Differential Equations and Dynamical Systems, coming to
Geometry through the seminal work of Andronov and  Pontrjagin \cite {ap} and
Peixoto \cite {mp}.

  The interest on lines of  harmonic mean curvature appears in   the
paper  of Raffy \cite {ra};
 see also  Eisenhart \cite {e}, section 55.
 The work of  Ogura \cite {og} regards these lines in terms   of his
unifying notion {\it T-Systems} and makes a local analysis of the
expressions of the fundamental quadratic forms in a chart whose coordinate
curves  are   lines of  harmonic mean curvature. A  comparative study of
these expressions with those corresponding to other lines  of geometric
interest, such as  the {\it principal, asymptotic, arithmetic}  and
{\it geometric mean curvature lines}   is carried out by Ogura in the
context of {\it T-Systems}, away from singularities.
 In the paper of  Occhipinti \cite {oc} is established the following
interesting projective relationship: {\it a line of harmonic mean
curvature divides harmonically those of geometric mean curvature (both) and
that (one) of arithmetic mean curvature }. See  \cite{ber}, chapter 6.

For being more descriptive and coherent with  that of  previous  recent
papers already cited, we adopt in this work  the denomination of {\it
harmonic mean curvature lines} instead of   {\it characteristic or
diagonal lines}, also  found  in the literature.

No global examples, or even local ones around singularities, of
harmonic mean curvature configurations seem to have been considered in the
literature on differential equations of classic differential geometry,
in contrast with the situations for the   principal and asymptotic cases
mentioned above. See also the work of  Anosov, for the global structure
of the geodesic flow \cite {an},   and that  of Banchoff, Gaffney and
McCrory \cite {bgm} for the parabolic and asymptotic lines.

\vskip 0.2cm

This paper is organized as follows:

Section 2 is devoted to the general study of the differential equations
and general properties of  Harmonic Mean Curvature Lines. Here are
given the precise definitions of the  Harmonic Mean Curvature Configuration
and of the  two transversal  Harmonic Mean Curvature Foliations with
singularities into which it splits.  The definition of  Harmonic Mean
Curvature Structural Stability focusing on the preservation  of the
qualitative properties of the foliations and the configuration under small
perturbations of the immersion, will be given at the end of this section.

In Section 3 the equation of lines of  harmonic mean curvature is
written in a Monge chart. The condition  for umbilic  harmonic mean
curvature stability  is explicitly stated in terms of the coefficients of the
third order jet of the    function which represents the immersion in a
Monge chart. The local  harmonic mean curvature separatrix
configurations at stable umbilics is established for $C^4$ immersions and resemble
the three  Darbouxian patterns of principal and arithmetic mean
curvature configurations \cite {da, gs1}. These patterns have been  also
recently established  for the case of  geometric mean curvature
configurations \cite {g}.

In Section 4 the derivative of first return Poincar\'e map along
a  harmonic mean curvature cycle is established. It consists of an
integral
expression
involving the curvature  functions
along the  cycle.

In Section 5 are studied the foliations by  lines of  harmonic
mean curvature  near the parabolic set   of an immersion, which
typically is   a regular curve. Three singular tangential patterns
exist generically in this case:  the {\it folded node}  the  {\it
folded saddle} and the {\it folded focus}. However, these types
alternate with the patterns established for the asymptotic lines
on the hyperbolic  region. The following is established and made
precise here: an elliptic harmonic (resp. asymptotic) saddle goes
adjacent with a hyperbolic (resp. harmonic) asymptotic node or
focus.   See subsection 5.1 and the pertinent bifurcation diagram.
Notice also that it has been proved that in the geometric mean
curvature  case   the {\it folded focus} is absent generically
\cite{g}.

Section 6 presents new  examples of  Harmonic Mean Curvature
Configurations on  the Torus of revolution  and  the quadratic Ellipsoid,
presenting non-trivial recurrences. This situation, impossible for   lines of
principal curvature, has been established, with different technical
details,  for arithmetic and geometric  mean curvature configurations in
\cite{m, g}.

In Section 7 the results presented in Sections 3, 4 and 5 are put
together
to provide sufficient conditions for  Harmonic Mean Curvature
Structural Stability. The density of these conditions is formulated and discussed at the
end of this section, however its rather technical proof will be postponed
to  another paper.

Section 8 contains an initial  discussion  motivated by this and
previous related papers.    We inquire about  the possibility and interest of
developing  a unifying general Theory for  Mean Curvature
Configurations,   valid for those already studied and also for possible  ``new" mean
curvature functions.

 \section{Differential Equations of   Harmonic Mean Curvature Lines}

Let $\alpha : {\mathbb M}^2\to \mathbb R^3$ be a $C^r,\;\; r\geq 4,$
immersion of
an oriented smooth surface ${\mathbb M}^2$ into $\mathbb R^3$. This
means that $D\alpha$
is injective  at every point in ${\mathbb M}^2$.

 The space $\mathbb R^3$ is oriented by a  once for all fixed
orientation
and endowed with the Euclidean inner product  $<,>$.

Let $N$  be a vector field orthonormal to $\alpha$.
Assume that $(u,v)$ is a positive chart of ${\mathbb M}^2$  and that
$\{\alpha_u, \alpha_v, N\}$
is  a positive frame in $\mathbb R^3$.

In the  chart $(u,v)$, the {\it first fundamental form} of an immersion
$\alpha$ is
given by:

$I_\alpha= <D\alpha,D\alpha>= E du^2+2Fdudv+Gdv^2$, with

$E=<\alpha_u,\alpha_u>$, $F=<\alpha_u,\alpha_v>$,
$G=<\alpha_v,\alpha_v>$

The  {\it second fundamental form} is given by:

 \centerline{$II_{\alpha}= <N , D^2\alpha> = e du^2 + 2f dudv + g
dv^2$.}

The normal curvature  at a point $p$ in a tangent direction $t=[du:dv]$
is given by:

 $$ k_n= k_n(p) =\frac{ II_\alpha(t,t)}{I_\alpha(t,t)}.$$

The lines of  harmonic   mean curvature  of $\alpha$ are regular curves
$\gamma$ on $\mathbb M^2$  having
normal  curvature   equal to the
 harmonic mean curvature of the immersion, i.e.,  $k_n=\frac{{\mathcal
K}}{{\mathcal H}}$, where $ {\mathcal K}$ = ${\mathcal K}{_\alpha}$ and
$ {\mathcal H}$ = ${\mathcal H}{_\alpha}$ are  the Gaussian and
Arithmetic Mean curvatures of $\alpha$.

Therefore the   pertinent differential equation for these lines is
given by:

$$\frac{edu^2+2fdudv+gdv^2}{Edu^2+2Fdudv+Gdv^2} =\frac{{\mathcal
K}}{{\mathcal H}}$$

Or equivalently by

\begin{equation}\label{eq:mgc}
 [g-\frac{{\mathcal K}}{{\mathcal H}}G]dv^2+2[f-\frac{{\mathcal
K}}{{\mathcal H}}F]dudv+  [e-\frac{\mathcal K}{{\mathcal
H}}E]du^2=0.\end{equation}
Also, as remarked by Occhipinti in   \cite{oc},  the equation of
harmonic curvature lines can be written as
 $$Jac(Jac(II,I),II)=0,$$ which leads to:

\begin{equation}\label{eq:harm}
 \aligned
L&dv^2  + M dudv + Ndu^2=0,\\
L=&  g(gE-eG) +2f(gF- fG) \\
M=&  2g(fE-eF) +2e(fG-gF)\\
N=&   e(eG-gE) +2f(f E-eF)    \endaligned
\end{equation}

This equation is defined only on the closure of the {\it Elliptic
region}, ${\mathbb E}{\mathbb M^2}{_\alpha}$,   of $\alpha$, where
${\mathcal K} >0$. It is bivalued and $C^{r-2},\;\; r\geq 4,$ smooth on the
complement of the umbilic,  ${\mathcal U}_\alpha$, and parabolic,
${\mathcal P}_\alpha$, sets of the immersion $\alpha$.  In fact, on ${\mathcal
U}_\alpha$ , where the principal curvatures  coincide, i.e where
${\mathcal H}^2-{\mathcal K}=0$, the equation vanishes identically; on
${\mathcal P}_\alpha$, it is univalued.

\vskip 0.3cm

 The developments above allow us
to organize the lines of  harmonic mean curvature of immersions  into
the
{\it   harmonic mean curvature
configuration,}  as follows:

 Through every point $p\in {\mathbb E}{\mathbb M^2}{_\alpha}\setminus
({\mathcal U}_\alpha \cup {\mathcal P}_\alpha)$,  pass two
 harmonic   mean curvature
lines of $\alpha$. Under the orientability
hypothesis imposed on $\mathbb M$, the  harmonic mean curvature
lines define two foliations:  $\mathbb H_{\alpha,1}$,
called the {\it minimal  harmonic mean curvature foliation}, along
which the
geodesic torsion is
negative (i.e  $\tau_g=- \sqrt{{\mathcal K}}\sqrt{  {\mathcal H}^2-
{\mathcal K} }/|{\mathcal H}|$ ),  and $\mathbb H_{\alpha,2}$, called the
{\it maximal
 harmonic mean curvature foliations}, along which the geodesic torsion
is
positive  (i.e  $\tau_g=\sqrt{{\mathcal K}}\sqrt{ {\mathcal H}^2-
{\mathcal K} }/|{\mathcal H}|$).

By comparison with the arithmetic mean curvature directions, making
angle $\pi/4$ with the minimal principal directions, the harmonic ones are
located  between them and the principal ones, making an angle
$\theta_h$ such that
$tan\theta_h$ =$\pm\sqrt{\frac{k_1}{k_2}}$,
 as follows from Euler's Formula. The particular expression for the
geodesic torsion given above
results from the formula  $\tau_g=(k_2-k_1)sin{\theta} cos{\theta}$
\cite{St},  is found in the work of Occhipinti \cite {oc}. See also Lemma
1 in Section 4 below. In \cite {oc, g} is also proved that  geometric
mean curvature lines are between the harmonic and arithmetic mean
curvature ones, making an angle $\theta_g$ such that $tan\theta_g$
=$\pm\sqrt{\frac{k_1}{k_2}}$.

With this data, Occhipinti \cite {oc}, has proved that
the two lines of mean geometric curvature, that of mean harmonic and
geometric curvature form a harmonic quadruple of lines.

The quadruple
$\mathbb H_\alpha=\{{\mathcal P}_\alpha,
{\mathcal U}_\alpha, \mathbb H_{\alpha,1},\mathbb H_{\alpha,2}\}$
is called the {\it   harmonic mean curvature configuration} of
$\alpha$.

It splits into two foliations with singularities:
$$\mathbb G^i_\alpha=\{{\mathcal P}_\alpha,
{\mathcal U}_\alpha, \mathbb H_{\alpha,i}\}, i=1, 2.$$

Let $\mathbb M^2$ be also   compact. Denote by ${\mathcal
M}^{r,s}({\mathbb M^2)}$ be the space of $ C^r $ immersions of $\mathbb M^2$ into
the Euclidean space $\mathbb R^3$, endowed with the $C^s$ topology.

 An immersion $\alpha$ is said $C^s$-{\it local   harmonic mean
curvature structurally
stable at a compact set} $C\subset \mathbb M^2$
if for any sequence of immersions $\alpha_n$ converging to $\alpha$ in
${\mathcal M}^{r,s}({\mathbb M^2)}$
there is a neighborhood $V_C$ of $C$, sequence of compact subsets $C_n$
and a sequence of homeomorphisms mapping $C$
to $C_n$ converging to the identity of $\mathbb M^2$ such that on $V_C$
it  maps umbilic and parabolic  points and  arcs of the  harmonic mean
curvature foliations
$\mathbb H_{\alpha,i}$  to those  of $\mathbb H_{\alpha_n,i}$
for $i = 1,\;2$.

An immersion $\alpha$ is said to be $C^s$-{\it   harmonic mean
curvature structurally stable} if the compact $C$ above is the closure of
${\mathbb E}{\mathbb M^2}{_\alpha}$.

Analogously, $\alpha$ is said  to be {\it i-} $C^s$-{\it   harmonic
mean curvature structurally stable}
if only the preservation of elements of  {\it i-th, i=1,2} foliation
with singularities is required.

A general study of the structural stability of  quadratic
differential equations (not necessarily derived from normal
curvature properties) has been carried out by Gu{\'\i}\~nez \cite
{gn}. See also the work of Bruce and Fidal  \cite {bf}
Bruce and Tari \cite{bt}, \cite{bt2} and Davydov \cite {dav}
 for the analysis of umbilic points  for general quadratic
 and also implicit
 differential equations.

For a study of the topology of  foliations with  non-orientable
singularities
on two dimensional manifolds,
see the works of Rosenberg and Levitt \cite{ro,le}. In these works the
leaves are not  defined by normal  curvature properties.

\section{  Harmonic mean curvature  lines near umbilic points }

 Let  $0$ be an umbilic point of a $ C^r ,\; r\geq  4,$ immersion
$\alpha$ parametrized in a Monge chart $(x,y)$ by
$\alpha(x,y)=(x,y,z(x,y))$, where

\begin{equation} \label{eq:1} h(x,y)=\frac k2(x^2+y^2)+ \frac a6
x^3+\frac b2 xy^2+\frac c6 y^3+O(4)\end{equation}

This reduced form is obtained by means of a rotation of the $x,y$-axes.
See \cite {gs1, gs2}.

According to Darboux \cite{da, gs1},  the differential equation of
principal curvature lines  is given by:

\begin{equation}
-[by+P_1]dy^2+[(b-a)x+cy+P_2]dxdy+[by+P_3]dx^2=0.\end{equation}

 As an starting point,   recall the behavior of principal lines near
Darbouxian umbilics
 in the following proposition.

\begin{proposition}{\rm \cite {gs1, gs2}} \label{prop:1} Assume  the
notation  established in \ref{eq:1}.
Suppose that the transversality condition $T:  b(b-a)\ne 0$  holds
 and
consider the following situations:
\begin{itemize}
\item[$D_1$)] $\;\;\Delta_{P}>0$

\item[$D_2$)] $\;\;   \Delta_{P} <0$ and $\dfrac ab >1$

\item[$D_3$)]   $\;\;    \dfrac ab  <1  $
\end{itemize}
 Here $\Delta_{P} =4b(a-2b)^3-c^2(a-2b)^2$

 Then each principal foliation  has  in a neighborhood of $0$,  one
hyperbolic sector in the $D_1$ case, one
parabolic and one hyperbolic sector in $D_2$ case and three hyperbolic
sectors in the case $D_3$.
  These  points  are called principal curvature  Darbouxian umbilics.
\end{proposition}

\begin{proposition}\label{prop:2}  Assume  the notation  established in
\ref{eq:1}.
Suppose that the transversality condition $T_h :  kb(b-a)\ne 0$  holds
 and
consider the following situations:
\begin{itemize}
\item[$H_1$)]   $\;\;\Delta_{h}>0$

\item[$H_2$)]$\;\;\Delta_{h}<0 $ and $ \;\;  \dfrac ab  >1$

\item[$H_3$)]    $\;\; \dfrac ab < 1. $
\end{itemize}
Here  $ \Delta_{h} =4c^2(2a-b)^2-[3c^2+(a-5b)^2][3(a-5b)(a-b)+c^2].  $

  {Then each   harmonic mean curvature foliation  has  in a
neighborhood of $0$,  one
hyperbolic sector in the $H_1$ case, one parabolic and one
hyperbolic sector in $H_2$ case and three hyperbolic sectors in
the case $H_3$. These umbilic points are called   harmonic mean
curvature Darbouxian umbilics.}

 The   harmonic mean curvature foliations $\mathbb H_{\alpha,i}$ near
an umbilic point
  of type $H_k$ has  a  local behavior  as shown in Figure 1.
  The {\it separatrices} of these singularities are called {\it umbilic
separatrices}.
\end{proposition}

 \begin{figure}[htbp] \label{fig:dgeo}
 \begin{center}
 \hskip 1cm
 \includegraphics[angle=0, width=11cm]{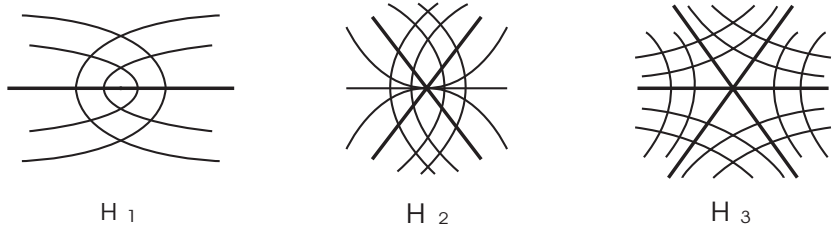}
 \caption{Harmonic mean curvature lines near the umbilic points $H_i$
and their separatrices}
   \end{center}
 \end{figure}

\begin{proof}
Near $0$, the functions $\mathcal K$ and $\mathcal H$ have the
following Taylor expansions.

$$  {\mathcal K}=k^2+  (a+b) k x+   ck y+O_1(2),\hskip .5cm
{\mathcal H}=k+ \frac 12(a+b)x+\frac 12 cy+O_2(2).$$

The  differential equation of the   harmonic mean curvature lines
\begin{equation} \label{eq:mh} [g-\frac{{\mathcal K}}{{\mathcal
H}}G]dv^2+2[f-\frac{{\mathcal K}}{{\mathcal H}}F]dudv+  [e-\frac{{\mathcal
K}}{{\mathcal H}}E]du^2=0 \end{equation}
\noindent  is given by:

\begin{equation}\label{eq:mg}
\aligned\; [(b-a)x+cy+M_1(x,y)]dy^2 +& [4by+M_2(x,y)]dxdy \\
\; -[(b-a)x  +  cy+M_3(x,y)]dx^2 =&0
\endaligned
\end{equation}

 \noindent where $M_i$,  $i=1,2,3$,  represent  functions of
order $O((x^2+y^2)).$

Thus, at the level of first jet, the differential equation \ref{eq:mg}
is the same as that of the arithmetic mean curvature lines given by
$$ [g- {\mathcal H}G]dv^2+2[f- {\mathcal H}F]dudv+  [e- {\mathcal
H}E]du^2=0,$$
\noindent as follows from the obvious fact that ${\mathcal H}$  and
$\frac{{\mathcal K}}{{\mathcal H}}$ have the same $1-$jet at $0$.

The conditions on $\Delta_h $ coincide with  those on   $\Delta_H$,
established to characterize the arithmetic mean curvature Darbouxian
umbilics studied in detail  in \cite{m}. Thus reducing the analysis of the
umibilic points to that of the  hyperbolicity  saddles and nodes whose
phase portrait is determined   only by  the first jet of the equation.
\end{proof}

\begin{theorem}\label{th:31}
An immersion $\alpha\in {\mathcal M}^{r,s}({\mathbb M^2)}$, $ r\geq
4$, is $C^3-$local
  harmonic mean curvature structurally  stable at ${\mathcal U}_\alpha$
if and only if
every $p\in {\mathcal U}_\alpha$ is one of the types $H_i$, $i=1,2,3$
of proposition \ref{prop:2}.
\end{theorem}

\begin{proof} Clearly proposition \ref{prop:2} shows that the condition
$H_i$, $i=1,2,3$ together with $T_h:kb(b-a)\ne 0 $  imply
the $C^3-$local   harmonic mean curvature structural stability. This
involves the construction of  the homeomorphism (by means of
canonical regions), mapping simultaneously  minimal and maximal
harmonic mean curvature lines  around the umbilic points of
$\alpha$ onto those of a $C^4$ slightly perturbed immersion.

We will discuss the necessity of the  condition $T_h:k(b-a)b\ne
0$ and of the conditions $H_i$, $i=1, \,2,\, 3$.  The first one
follows from
its identification with a transversality condition that guarantees the
persistent isolatedness  of the umbilic
points of $\alpha$ and its separation from the parabolic set, as well
as the persistent regularity of  the Lie-Cartan
  surface
${\mathcal G}$, obtained from the projectivization of the equation
\ref{eq:mh}.  Failure of $T_h$
condition has the following implications:
\begin{itemize}
\item[a)] $b(b-a)=0$; in this case the elimination or splitting of the
umbilic point can be achieved by small perturbations.

\item[b)] $k=0$ and $b(b-a)\ne 0$; in this case a small perturbation
separates the umbilic point from the parabolic set.
\end{itemize}
 The necessity of condition  $H_i$ follows from its dynamic
identification with the hyperbolicity  of the
equilibria along the projective line of the vector field obtained
lifting equation
 (\ref{eq:mh}) to the surface   ${\mathcal G}$. Failure of this
condition would
 make possible
to  change the number of  harmonic mean curvature umbilic separatrices
at the umbilic point by means a small  perturbation of
the immersion.
\end{proof}

\section{Periodic   Harmonic Mean Curvature Lines}

Let $\alpha:\mathbb M^2\to \mathbb R^3$ be an immersion of a compact
and
 oriented surface and consider the foliations $\mathbb H_{\alpha, i}$,
$i=1,\;2$, given by
 the {\it   harmonic mean curvature lines}.

  In terms of geometric invariants, here is established an integral
expression  for  the first derivative of the  return map of a
periodic    harmonic mean curvature line, called {\it   harmonic mean
curvature cycle}.
Recall that the return map associated  to a cycle
is a local diffeomorphism with a fixed point, defined on a cross
section  normal to the cycle by following  the
integral curves through this section until they meet again the section.
This map is called holonomy in   Foliation
Theory and Poincar\'e Map in Dynamical Systems, \cite{pm}.

A   harmonic mean curvature cycle is called {\it hyperbolic} if the
first derivative of the return map at the fixed  point is
different from one.

{ The   harmonic mean curvature foliations $\mathbb H_{\alpha,i}$ has
no  harmonic mean curvature cycles such that  the return map reverses
the orientation.}
 Initially, the integral expression for the derivative of the return
map is
obtained in class $ C^6$; see Lemma \ref{lm:42} and Proposition
\ref{prop:41}. Later on,
in Remark
\ref{rm:41} it is shown how to extend it to class $C^3$.

The characterization of hyperbolicity of  harmonic mean curvature
cycles in terms
of local structural stability is given in  Theorem  \ref{th:41} of
this section.

\begin{lemma}\label{lm:41} Let $c:I\to {\mathbb M}^2$ be a  harmonic
mean curvature
line parametrized by arc length. Then the Darboux frame is given by:

$$\aligned
 T^\prime &= k_g N\wedge T+ \frac{{\mathcal K}}{{\mathcal H}} N\\
(N\wedge T)^\prime &= -k_g T+ \tau_g N\\
N^\prime &= -\frac{{\mathcal K}}{{\mathcal H}} T-\tau_g N\wedge
T\endaligned\eqno
$$
 \noindent where $\tau_g=  \pm \sqrt{  {\mathcal K}} \frac{ \sqrt{
{\mathcal H}^2 - {\mathcal K}}}{|{\mathcal H}|}  $. The sign of $\tau_g$
is positive (resp. negative) if $c$ is maximal (resp. minimal)
harmonic mean curvature line.
\end{lemma}

\begin{proof} The normal curvature $k_n$ of the curve $c$ is by the
definition
 the   harmonic mean curvature $\frac{\mathcal K}{\mathcal H}$.
From the Euler equation
$k_n=k_1\cos^2\theta+k_2\sin^2\theta=\frac{\mathcal K}{{\mathcal H}}$, get
$\tan\theta=\pm\sqrt{\frac{k_1}{k_2}}$.    Therefore, by direct
calculation, the geodesic torsion is given by
$\tau_g=(k_2-k_1)\sin\theta\cos\theta  =  \pm \sqrt{  {\mathcal K}}
\frac{\sqrt{{\mathcal H}^2 - {\mathcal K}}}{\mathcal H}$.
\end{proof}

 \begin{remark}
The expression for the geodesic curvature $k_g$ will not be needed
explicitly in this work. However, it can be given in terms of the principal
curvatures and their derivatives using a formula due to Liouville
\cite{St}.
\end{remark}

\begin{lemma}\label{lm:42} Let $\alpha:\mathbb M\to \mathbb R^3$ be an
immersion of class $ C^r $, $ r\geq   6$, and $c$ be a mean
curvature cycle of $\alpha$,
 parametrized by arc length and  of length $L$. Then the expression,

$$\alpha(s,v)=c(s)+ v(N\wedge T)(s)+[ (2{\mathcal H} (s)-
\frac{{\mathcal K}}{{\mathcal H}}(s))   \frac{v^2}2+
 \frac{A(s)}{6}v^3+ v^3B(s,v)]N(s)\eqno$$

 \noindent where $B(s,0)=0$, defines    a local chart $(s,v)$ of class
$C ^{r-5}$ in a neighborhood of $c$.
\end{lemma}

\begin{proof}  The curve $c$ is of class $C ^{r-1}$ and the map
$\alpha(s,v,w)=c(s)+ v(N\wedge T)(s)+wN(s)$ is of class $C ^{r-2}$ and
is a local diffeomorphism  in a neighborhood of the axis $s$. In fact
$[\alpha_s,\alpha_v,\alpha_w](s,0,0)=1$. Therefore there is a function
$W(s,v)$ of class $C ^{r-2}$ such that $\alpha(s,v,W(s,v))$ is a
parametrization of a tubular neighborhood of $\alpha\circ c$. Now for each
$s$, $W(s,v)$ is just a parametrization of
the curve of intersection between $\alpha(\mathbb M)$ and  the normal
plane generated by $\{(N\wedge T)(s), N(s)\}$. This curve of
intersection  is tangent to $(N\wedge T)(s)$ at $v=0$ and notice that
$k_n(N\wedge T)(s)=2{\mathcal H}(s)-\frac{{\mathcal K}}{{\mathcal H}}(s) $.
Therefore,

\begin{equation} \aligned  \alpha(s,v,W(s,v))=& c(s)+ v(N\wedge T)(s)\\
+&[(2{\mathcal H} (s)- \frac{{\mathcal K}}{{\mathcal H}}(s))\frac{
v^2}2+
 \frac{A(s)}6v^3+ v^3B(s,v)]N(s), \endaligned \end{equation}

 \noindent where $A$ is of class $C ^{r-5}$ and $B(s,0)=0$.
\end{proof}

We now  compute the coefficients of the first and second fundamental
forms in  the chart $(s,v)$ constructed above, to be used in  proposition
\ref{prop:41}.

$$\aligned
N(s,v)=&\frac{\alpha_s\wedge\alpha_v}{\mid\alpha_s\wedge\alpha_v\mid}
 = [-\tau_g(s)v+O(2)]T(s)\\
-&[(2{\mathcal H} (s)- \frac{{\mathcal K}}{{\mathcal
H}}(s))v+O(2)](N\wedge T)(s)
+[1+ O(2)]N(s).\endaligned $$

 \noindent Therefore it follows that  $E=<\alpha_s,\alpha_s>$,
$F=<\alpha_s,\alpha_v>$, $G=<\alpha_v,
\alpha_v>$,
$e=<N,\alpha_{ss}>,$
 $\;\; f=<N,\alpha_{sv}>\;\;$ and $\; g=<N,\alpha_{vv}>$ are given by

\begin{equation}\label{eq:1f2f}
\aligned
E(s,v) &= 1-2k_g(s)v+h.o.t\\
F(s,v) &=  0+ 0.v+h.o.t\\
G(s,v) &= 1 +0.v+h.o.t\\
e(s,v)&=\frac{{\mathcal K}}{{\mathcal H}}(s)
+v[\tau_g^\prime(s)-2k_g(s) {\mathcal H}(s) ]+ h.o.t\\
f(s,v) &= \tau_g(s)+ \{[2{\mathcal H}(s)-\frac{{\mathcal K}}{{\mathcal
H}}(s)]^\prime +k_g(s)\tau_g(s)\}v+ h.o.t\\
g(s,v) &=2{\mathcal H}(s)-\frac{{\mathcal K}}{{\mathcal H}}(s) +
A(s)v+ h.o.t \endaligned\end{equation}

\begin{proposition}\label{prop:41} Let $\alpha:\mathbb M\to \mathbb
R^3$ be an immersion of class $ C^r $, $ r\geq   6$ and $c$ be closed
harmonic line $c$  of $\alpha$,
 parametrized by arc length $s$ and  of total length $L$.
 Then the derivative of the Poincar\'e map $\pi_\alpha$ associated to
$c$ is given by:

$$ln\pi_{\alpha}^\prime(0)=  \int_0^L\left[{\frac{[\frac{\mathcal
K}{{\mathcal H}}]_v  }{2\tau_g}} +\frac{ k_g}{\tau_g}({\mathcal H}-
\frac{\mathcal K}{\mathcal H})   \right]ds. \eqno$$

\noindent Here $\tau_g$=$\pm{\frac{\sqrt{\mathcal K}}{\mathcal
H}}\sqrt{{\mathcal H}^2-{\mathcal K}}$.

\end{proposition}

\begin{proof} The Poincar\'e map associated to $c$ is the map
$\pi_\alpha:\Sigma \to \Sigma$ defined in a transversal section to $c$ such that
$\pi_{\alpha}(p)=p$ for $p\in c\cap\Sigma$ and $\pi_{\alpha}(q)$ is the
first return of the  harmonic mean curvature line through $q$ to the
section $\Sigma$, choosing a positive orientation for $c$. It is a local
diffeomorphism and is defined, in the local chart $(s,v)$ introduced in
Lemma \ref{lm:42}, by $\pi_\alpha:\{s=0\}\to \{s=L\}$,
$\pi_\alpha(v_0)=v(L,v_0)$, where $v(s,v_0)$ is the solution of the Cauchy problem

$$(g-\frac{{\mathcal K}}{{\mathcal H}})dv^2+2(f-\frac{{\mathcal
K}}{{\mathcal H}} F)dsdv+(e-\frac{{\mathcal K}}{{\mathcal H}}E)ds^2=0, \quad
v(0,v_0)=v_0.$$

Direct calculation gives
that the derivative  of the Poincar\'e map satisfies
 the following linear differential equation:

$$\frac{d}{ds}(\frac{dv}{dv_0})
 = -\frac{N_v}M(\frac{dv}{dv_0})=-\frac{ [ e-\frac{{\mathcal
K}}{{\mathcal H}}(s) E] _v}{  2[ f-\frac{{\mathcal K}}{{\mathcal H}}(s) F]  }
  (\frac{dv}{dv_0})\eqno$$

Therefore, using equation \ref{eq:1f2f} it results that
$$\frac{ [ e-\frac{{\mathcal K}}{{\mathcal H}}(s) E] _v}{  2[
f-\frac{{\mathcal K}}{{\mathcal H}}(s) F]  } =
-\frac{\tau_g^\prime}{2\tau_g}-\frac{[\frac{{\mathcal K}}{{\mathcal H}}(s) ]_v  }{2\tau_g}
-\frac{k_g}{\tau_g}({\mathcal H}- \frac{\mathcal K}{{\mathcal H}}). $$

Integrating  the equation above  along an arc
$[s_0,s_1]$ of  harmonic mean
curvature line,  it follows that:

 \begin{equation}\label{eq:47}
\frac{dv}{dv_0}|_{v_0=0}= \frac
{(\tau_g(s_1))^{\frac{-1}2}}{(\tau_g(s_0))^{\frac{-1}2}}exp[
\int_{s_0}^{s_1} \left[\frac{[\frac{\mathcal K}{{\mathcal H}}]_v
}{2\tau_g} +\frac{ k_g}{\tau_g}({\mathcal H}- \frac{\mathcal K}{{\mathcal
H}})   \right]ds.\end{equation}

Applying \ref{eq:47} along the  harmonic mean curvature cycle of length
$ L$, obtain
 $$\frac{dv}{dv_0}|_{v_0=0}=  exp[
  \int_0^L\left[\frac{[\frac{\mathcal K}{{\mathcal H}}]_v  }{2\tau_g}
+\frac{ k_g}{\tau_g}({\mathcal H}- \frac{\mathcal K}{{\mathcal H}})
\right]ds. $$

From the equation ${\mathcal K}=(eg-f^2)/(EG-F^2)$ evaluated
at $v=0$  it follows that ${\mathcal K}=\frac{\mathcal K}{{\mathcal
H}}[2{\mathcal H}-\frac{\mathcal K}{{\mathcal H}}]-\tau_g^2.$
Solving this equation it follows that
$\tau_g$=$\pm{\frac{\sqrt{\mathcal K}}{\mathcal H}}\sqrt{{\mathcal
H}^2-{\mathcal K}}$

 \noindent This ends the proof. \end{proof}

  \begin{remark}\label{rm:41}  At this point we show how to extend
 the expression for the derivative of the hyperbolicity of
 harmonic mean curvature cycles established for class  $C^6$ to  class
$C^3$
 (in fact we need only class $C^4$).

The expression \ref{eq:47}
 is the derivative of the transition map
for a  harmonic mean curvature foliation (which at this point is only
of class
$C^1$), along an arc of  harmonic mean curvature line. In fact, this
follows by
approximating  the $C^3$ immersion by one of class  $C^6$. The
corresponding
transition map (now of class $C^4$) whose derivative is given by
expression \ref{eq:47}
converges to the original one (in
class $C^1$) whose expression must given by the same integral, since
the functions involved there are the uniform limits of the
corresponding
ones for the approximating immersion.
\end{remark}

\begin{proposition}\label{prop:42} Let $\alpha:\mathbb M\to \mathbb
R^3$ be an immersion of class $ C^r $, $ r\geq   6$,  and  $c$ be a
maximal   harmonic mean curvature cycle of $\alpha$, parametrized by
arc length
and of length $L$. Consider a chart $(s,v)$ as in lemma \ref{lm:42} and
consider the deformation
$$\beta_\ep(s,v)=\beta(\ep,s,v)=\alpha(s,v)+\ep
[\frac{A_1(s)}6v^3]\delta (v)N(s)\eqno$$
 \noindent where $\delta =1$ in neighborhood of $v=0$, with small
support and $A_1(s)=\tau_g(s)>0$.

Then $c$ is a  harmonic mean curvature cycle of $\beta_\ep$ for all
$\ep $ small and $c$ is a hyperbolic
 harmonic mean curvature cycle for $\beta_\ep$, $\ep\ne 0$.
\end{proposition}

\begin{proof}  In the chart $(s,v)$, for the immersion $\beta_\ep$, it
is obtained that:

$$\aligned
E_\ep(s,v) &= 1-2k_g(s)v+h.o.t\\
F_\ep(s,v) &=  0+ 0.v+h.o.t\\
G_\ep(s,v) &= 1 +0.v+h.o.t\\
e_\ep(s,v)&=\frac{{\mathcal K}}{{\mathcal H}}(s)
+v[\tau_g^\prime(s)-2k_g(s){\mathcal H}(s)\; )]+ h.o.t\\
f_\ep(s,v) &= \tau_g(s)+ [(2{\mathcal H}(s)-\frac{{\mathcal
K}}{{\mathcal H}}(s))^\prime +k_g\tau_g]v+ h.o.t\\
g_\ep(s,v) &=2{\mathcal H}(s)-\frac{{\mathcal K}}{{\mathcal H}} (s)
+v[A(s)+\epsilon A_1(s)]+ h.o.t \endaligned\eqno$$

In the expressions above $E_\ep=<\beta_s,\beta_s>$,
$F_\ep=<\beta_s,\beta_v>$, $G_\ep=<\beta_v,\beta_v>$,
 \noindent $\; e_\ep=<\beta_{ss},N>$,
 $\; f_\ep=<N,\beta_{sv}>,\;$  $\; g_\ep=<N,\beta_{vv}>$, where
$N=N_\ep=\beta_s\wedge \beta_v /\mid\beta_s\wedge\beta_v\mid.$

   For all $\ep$ small it follows  that:

$$\aligned (e_\ep -\frac{{\mathcal K}_\ep}{{\mathcal H}_\ep}
E_\ep)(s,0,\epsilon)=& 0\\
{{\mathcal K}_\ep}_v(s,0,\epsilon)=&  \ep \frac{{\mathcal
K}_\ep}{{\mathcal H}_\ep}A_1(s)+ f_1(k_g, \tau_g,  {\mathcal K} ,  {\mathcal
H})(s)\\
{{\mathcal H}_\ep}_v(s,0,\epsilon)=& \frac 12  \ep  A_1(s)+ f_2(k_g,
\tau_g,  {\mathcal K} ,  {\mathcal H})(s)\\
\frac{d}{d\ep }\big[\frac{{\mathcal K}_\ep}{{\mathcal
H}_\ep}\big]_v|_{\epsilon=0}=&\frac 12  \frac{{\mathcal K}}{{\mathcal H}^2} A_1(s).
\endaligned
$$
Therefore  $c$ is a maximal   harmonic mean curvature cycle for all
$\beta_\ep$.

Assuming that  $A_1(s)=4\tau_g(s)>0$,  it results that
$$\frac{d}{d\ep}(ln\pi^\prime(0))|_{\ep=0}=  \int_0^L
\frac{d}{d\ep}\left(\frac{(\frac{{\mathcal K}_\ep}{{\mathcal
H}_\ep})_v}{2\tau_g}+\frac{ k_g}{\tau_g}({{\mathcal H}_\ep}- \frac{{\mathcal
K}_\ep}{{\mathcal H}_\ep})\right)ds=\int_0^L \frac{{\mathcal K}}{{\mathcal
H}^2}ds   > 0.  $$
\end{proof}

As a synthesis of propositions \ref{prop:41} and \ref{prop:42},  the
following theorem is obtained.

\begin{theorem}\label{th:41}
An immersion $\alpha\in {\mathcal M}^{r,s}({\mathbb M^2)}$, $ r\geq
6$, is  $C^6-$local  harmonic mean curvature structurally stable at a
harmonic mean curvature cycle $c$
if only if,
$$  \int_0^L\left[\frac{[\frac{\mathcal K}{{\mathcal H}}]_v  }{2\tau_g}
+\frac{ k_g}{\tau_g}({\mathcal H}- \frac{\mathcal K}{{\mathcal H}})
\right]ds \neq 0. $$
\end{theorem}

\begin{proof} Using propositions  \ref{prop:41} and \ref{prop:42}, the
local topological character of the foliation can be changed by small
perturbation of the immersion, when the cycle is not hyperbolic.
\end{proof}

\section{ Harmonic Mean Curvature Lines near the Parabolic Curve}

 Let  $0$ be a parabolic  point of a $ C^r , \; r\geq   6$, immersion
$\alpha$ parametrized in a Monge chart $(x,y)$ by
$\alpha(x,y)=(x,y,z(x,y))$, where

\begin{equation}\aligned z(x,y) =&\frac k2 y^2 + \frac a6 x^3+\frac b2
xy^2+\frac d2 x^2y+\frac c6 y^3\\ +&\frac A{24}x^4+\frac B{6}
x^3y+\frac{C}4 x^2y^2+\frac D6 xy^3+\frac E{24}y^4+O(5)\endaligned
\end{equation}

The coefficients of the first and second fundamental forms are given
by:

\begin{equation}\label{eq:1f2fp}\aligned E(x,y)=& 1+O(4)\\
F(x,y)=&  +O(3)\\
G(x,y)=&1+k^2 y^2 +O(3)\\
e(x,y)=& ax+dy+\frac A2x^2+Bxy+\frac C2y^2 +O(3)\\
f(x,y)=& dx+ by+\frac B2x^2+Cxy+\frac D2 y^2 +O(3)\\
g(x,y)=& k+bx+cy+\frac C2x^2+Dxy+\frac 12{( E-k^3)} y^2
+O(3)\endaligned
\end{equation}

The Gaussian   and the Arithmetic Mean curvatures are given by

\begin{equation}\label{eq:gc}
 \aligned
{\mathcal K}(x,y)=& k(ax+dy)+ \frac 12(Ak+2ab-2d^2)x^2+(Bk+ac-bd)xy\\
+&\frac 12(Ck+2cd-2b^2)y^2+O(3),\\
{\mathcal H}(x,y)=& \frac 12k+\frac 12(a+b)x+\frac 12 (c+d)y+(A+C
)\frac{x^2}4\\
+& (B+D)xy +( E -3 k^3+ C) \frac{y^2}4+O(3)
\endaligned \end{equation}

The coefficients of the quadratic differential equation \ref{eq:harm}
are given by

\begin{equation}\label{eq:ai}
\aligned
L=& k^2+k(2b-  a)x+k( 2c-  d) y +(2kC - Ak+  2b^2+4d^2-2 ab)\frac{x^2}2
\\
+& ( 3 db - ac+2kD- kB+2cb) xy \\
+&( 2 c^2+4b^2+ 2kE- 2 cd-  kC-2k^4)\frac{y^2}2 +O(3)\\
M =& 2k (d.x+b.y) +(4ad+2kB+4bd)\frac{x^2}2+2(b^2+d^2+ab+kC+cd)xy\\
+&   ( 4bd+ 2kD+4cb ) \frac{y^2}2+O(3)\\
 N=&- k(a x + d  y) +(2 a^2+ 4d^2- 2 ab- A k)\frac{x^2}2 \\
+& (2a d - kB+ 3bd - ac) xy +(2 d^2+4b^2- kC-2c d )\frac{y^2}2+O(3)\\
\endaligned
\end{equation}

\begin{lemma}\label{lm:p1} Let $0$ be a parabolic point and consider
the parametrization $(x,y,h(x,y))$ as  above. If $k>0$ and $a^2 + d^2 \ne
0$ then the set of parabolic points is locally a regular curve normal
to the vector $(a,d)$ at $0$.

 If $a\ne 0$ the parabolic curve is transversal to the minimal
principal direction $(1,0)$.

If $a=0$ then the parabolic curve is tangent to the principal direction
given by $(1,0)$ and has quadratic contact with the corresponding
minimal principal curvature line if $dk(Ak-3d^2)\ne 0$.

\end{lemma}

\begin{proof} If $a \ne 0$, from the expression of $\mathcal K$ given
by equation \ref{eq:gc} it follows that the parabolic line is given by
$x=-\frac{d}{a}y+O_1(2)$ and so is transversal to the principal
direction $(1,0)$ at $(0,0)$.

 If $a=0$, from the expression of $\mathcal K$ given by equation
\ref{eq:gc} it follows that the parabolic line is given by $y=\frac{
2d^2-Ak}{2dk}x^2+O_2(3)$  and that $y=-\frac{d}{2k}x^2+O_3(3)$ is the principal
line tangent to the principal direction $(1,0)$.
Now the condition of quadratic contact $\frac{ 2d^2-Ak}{2dk}\ne
-\frac{d}{2k}$ is equivalent to $dk(Ak-3d^2)\ne 0$.\end{proof}

\begin{proposition}\label{prop:oc-p} Let $0$ be a parabolic point and
the Monge chart $(x,y)$ as above.

If $a\ne 0$ then the mean   harmonic curvature lines are transversal to
the parabolic curve and the mean curvatures lines are shown in the
picture below, the cuspidal case.

If $a=0$ and $\sigma=  k^2(A k-3 d ^2)\ne 0$ then the  mean   harmonic
curvature lines are shown in the picture below. In fact, if
$\sigma >0$ then the mean   harmonic curvature lines are folded
saddles. Otherwise,  if $\sigma < 0$ then the mean   harmonic
curvature lines are folded nodes or folded focus according to
$\delta=-23 d ^2+8 A k  $ be positive or negative. The two separatrices of these
tangential singularities, folded saddle and folded node, as
illustrated in the Figure 2 below, are called parabolic
separatrices.
\end{proposition}

 \begin{figure}[htbp] \label{fig:hparab}
 \begin{center}
 \includegraphics[angle=0, width=11cm]{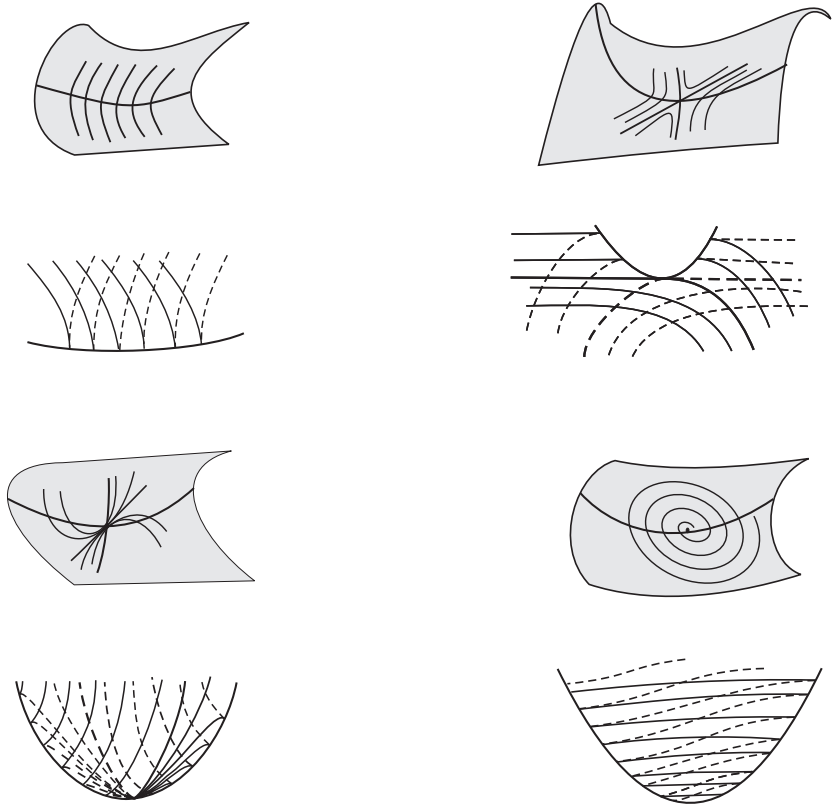}
 \caption{Harmonic mean curvature lines near a parabolic point
(cuspidal, folded saddle,    folded node and folded focus) and their
separatrices}
    \end{center}
  \end{figure}

\begin{proof}
Consider the quadratic differential equation
$$H(x,y,[dx:dy])=Ldy^2+Mdxdy +Ndx^2 =0$$

\noindent and the Lie-Cartan line field $X$ of class $C ^{r-3}$ defined
by

$$\aligned x^\prime =&H_p\\
y^\prime =& pH_p\\
p^\prime =&-(H_x+pH_y), \hskip 1cm p=\frac{dy}{dx}\endaligned$$
\noindent where $L$, $M$ and $N$ are given by equation \ref{eq:ai}.

If $a\ne 0$ the vector $Y$ is regular and therefore  the mean
harmonic curvature lines are transversal to the parabolic line and at
parabolic points these lines are  tangent to the  principal direction $(1,0)$.

If $a=0$, direct calculation  gives $H(0)=0,\;\; H_x(0)=0, \;\;
H_y(0)=-kd, \;\; H_p(0)=0.$

 \begin{equation}\aligned
DX(0)=\left(\begin{matrix} 2kd & 2kb & 2k^2\\  0 & 0 & 0\\    Ak-4d^2 &
kB- 3 bd  & -  kd\end{matrix}\right)
\endaligned
\end{equation}

The non vanishing eigenvalues of $DX(0)$ are
$$\lambda_1=(\frac 12 d+\frac 12\sqrt{-23d^2+8Ak})k,\;\;
\lambda_2=(\frac 12 d-\frac 12\sqrt{-23d^2+8Ak})k$$
Therefore, $\lambda_1\lambda_2=- 2k^2(Ak-3d^2)$.

It follows that $0$ is a hyperbolic  singularity provided $\sigma
(Ak-3d^2)kd\ne 0$. If $\sigma >0$ then the mean   harmonic curvature lines
are folded saddles and   if $\sigma <0$ then the   mean   harmonic
curvature lines are folded nodes $(8Ak-23d^2>0 )$ or folded focus
$(8Ak-23d^2<0)$.   See Figure 2 above.

 \end{proof}

\begin{theorem}\label{th:51}
An immersion $\alpha\in {\mathcal M}^{r,s}({\mathbb M^2)}$, $ r\geq
6$, is  $C^6-$local  harmonic mean curvature structurally stable at a
tangential parabolic point $p$
if only if, the condition $\sigma\delta \neq 0$ in  proposition
\ref{prop:oc-p}  holds.

\end{theorem}

\begin{proof}
Direct from Lemma   \ref{lm:p1} and proposition \ref{prop:oc-p}, the
local topological character of the foliation can be changed by small
perturbation of the immersion when $\delta \sigma = 0$.
\end{proof}

\subsection{ Asymptotic Lines near a Parabolic Curve}

\begin{proposition}\label{prop:oc-a} Let $0$ be a parabolic point and
the Monge chart $(x,y)$ as above.

If $a\ne 0$ then the mean   asymptotic lines are transversal to the
parabolic curve and   are shown in the picture below, the cuspidal case.

If $a=0$ and $\sigma=  k^2(A k-3 d ^2)\ne 0$ then the asymptotic are
shown in the picture below. In fact, if
$\sigma <0$ then the asymptotic  lines are folded
saddles. Otherwise,  if $\sigma >0$ then the asymptotic lines are
folded nodes or folded focus according to $\delta_a=25 d ^2-8 A k  $ be
positive or negative. The two separatrices of these
tangential singularities, folded saddle and folded node, as
illustrated in  Figure 3 below, are called parabolic
separatrices.
\end{proposition}

\begin{figure}[htbp] \label{fig:ha}
 \begin{center}
 \includegraphics[angle=0, width=9cm]{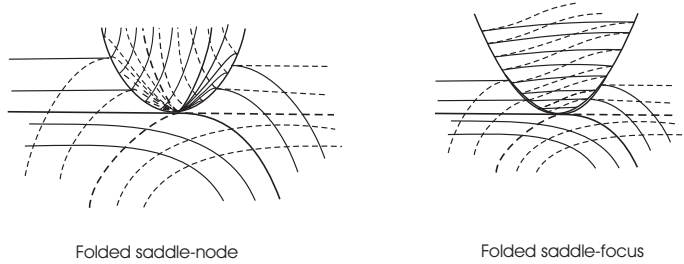}
 \caption{Harmonic Node and Focus adjacent to Asymptotic Saddle.}
   \end{center}
 \end{figure}

\begin{proof} The proof follows from direct calculations similar to
those performed in proposition \ref {prop:oc-p}.
In fact, considering the implicit differential equation
$${\mathcal A}(x,y,p)= gp^2+2fp+e=0, \;\;\;\; p=\frac{dy}{dx}$$
\noindent where $e$, $f$ and $g$ are given by equation \ref{eq:1f2fp}
and the Lie-Cartan line field
$$Y={\mathcal A}_p\frac{\partial }{\partial x}+p{\mathcal
A}_p\frac{\partial }{\partial y}-({\mathcal A}_x+p{\mathcal A}_y)\frac{\partial
}{\partial x},$$
\noindent it follows that

 \begin{equation}\aligned
DY(0)=\left(\begin{matrix} 2 d & 2 b & 2k \\  0 & 0 & 0\\    -A & -B  &
-  3d\end{matrix}\right)
\endaligned
\end{equation}
The non vanishing eigenvalues of $DY(0)$ are
$$r_1= \frac 12 d+\frac 12\sqrt{ 25d^2-8Ak} ,\;\; r_2= \frac 12 d-\frac
12\sqrt{ 25d^2-8Ak} $$
Therefore, $r_1 r_2= 2(Ak-3d^2)$.

It follows that $0$ is a hyperbolic  singularity provided $  Ak-3d^2
\ne 0$. If $Ak-3d^2  <0$ then the mean   harmonic curvature lines are
folded saddles;  if $Ak-3d^2  >0$ then the   mean   harmonic curvature
lines are folded nodes $( 25d^2-8Ak>0 )$ or folded focus
$(25d^2-8Ak<0)$.   See Figure 3 above.
 \end{proof}

\begin{remark} The geometric conditions of  asymptotic folded saddles,
nodes and focus near a parabolic line  was obtained in \cite{a1}.
\end{remark}

\begin{remark} In the plane $k=1$ the diagram of folded saddles, folded
nodes and folded focus for harmonic mean curvature lines and asymptotic
lines is as shown in  Figure 4 below.
\end{remark}

\begin{figure}[htbp] \label{fig:dbif}
 \begin{center}
 \includegraphics[angle=0, width=12cm]{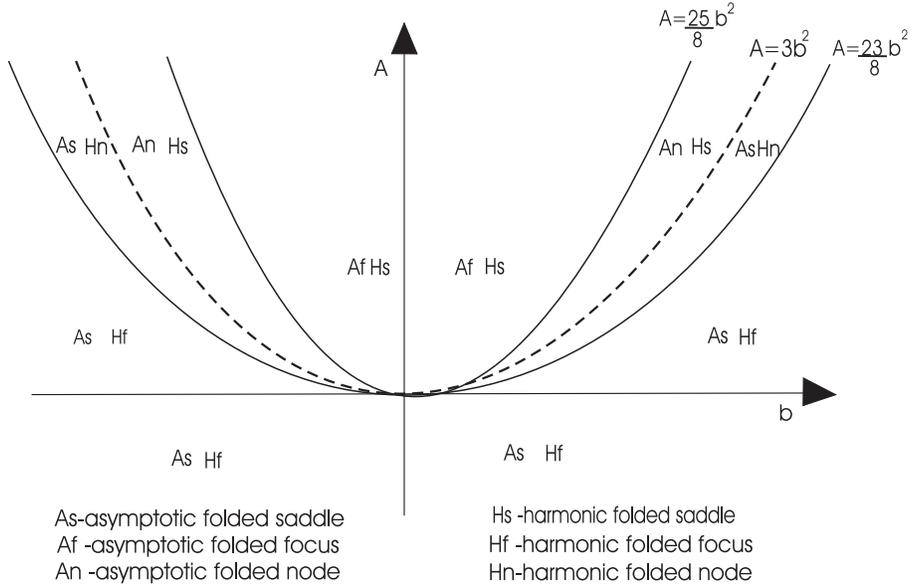}
 \caption{Bifurcation diagram of asymptotic
 and harmonic mean curvature lines in the plane b$\times$A}
   \end{center}
 \end{figure}

\section{Examples of Harmonic Mean Curvature Configurations}

As mentioned in the Introduction,  no examples of  harmonic mean
curvature foliations are given in the literature, in contrast with the
principal and asymptotic foliation. In this section are studied the  harmonic
mean curvature configurations in two classical surfaces: The Torus and
the Ellipsoid.
In contrast with the principal case \cite{Sp, St} (but in concordance
with the arithmetic mean curvature one \cite{m}) non-trivial recurrence
can occur here.

\begin{proposition}\label{prop:mgtoro}
Consider a torus of revolution $T(r,R)$ obtained by rotating a
circle of radius $r$ around a line in the same plane and at a
distance $R$, $R>r$, from its center.
 Define the function $\rho$ of $a=\frac rR$, as follows:
$$\rho= \rho(a)=\int_{-\frac{\pi}2}^{\frac{\pi}2}\sqrt{\frac{a} {\cos
s(1+a\cos s)}}ds.$$
  Consider the regular curves (folded extended harmonic lines) defined as
  the union  of harmonic lines and parabolic points ( a harmonic line of one
foliation that arrive  at the parabolic set at a given  point
is continued through the line  of the other foliation leaving the
parabolic set at  this point and so on). Then the  folded extended
 harmonic mean curvature lines on $T(r,R)$, defined in
the elliptic region are all closed or all recurrent according  to
$\rho\in \mathbb Q$ or $\rho \in \mathbb R\setminus\Bbb Q$.
Furthermore, both cases occur for appropriate $(r,R)$.
\end{proposition}

\begin{proof} The torus of revolution $T(r,R)$ is parametrized by

 $$\alpha(s,\theta)=((R+r\cos s )\cos \theta ,(R+r\cos s
)\sin \theta, \; r\sin s ).$$

 \noindent Direct calculation shows that $E=r^2,\;$ $\;F=0$,
$\;G=[R+r\cos s ]^2,\;$  $e=-r$, $\;f=0$ and
 $\;g=-\cos s(R+r\cos s)$.
  Clearly  $(s,\theta)$ is a principal chart.

 The differential equation of the  harmonic mean curvature lines, in
the principal chart $(s,\theta)$, is given by
$ e ds^2-g d\theta^2=0$.
This is equivalent to
$$-\cos s (1+a\cos s) d\theta^2+ a ds^2=0, \;\;\;\; a=\frac rR$$

Solving the equation above it  follows that,

$$\int_{\theta_0}^{\theta_1} d\theta= \pm\int_{s_0}^{s_1}\sqrt{\frac{a}
{\cos s(1+a\cos s)}}ds
   .$$ So the  two Poincar\'e maps,
$\pi_\pm:\{s=-\frac{\pi}2\}\to
\{s=\frac{\pi}2\}$, defined by $\pi_\pm(\theta_0)=\theta_0 \pm
 2\pi\rho(\frac rR)$ have rotation number equal to $\pm \rho(\frac rR)$.
Direct calculations gives that $\rho(a)\ne  0$ and $\rho^\prime(a)\ne
0$ for $a>0$.
Therefore,  both the rational and irrational cases occur. This
ends the proof.
\end{proof}

\begin{proposition}\label{prop:mgelipsoide}   Consider the  ellipsoid
$\mathbb E_{a,b,c}$
with three axes $a>b>c>0$. Then $\mathbb E_{a,b,c}$ have four
umbilic points located in the plane of symmetry orthogonal to
middle axis; they are of the type $H_1$ for  harmonic mean
curvature lines and of type $D_1$ for the principal curvature
lines.\end{proposition}

\begin{proof} This follows from proposition \ref{prop:2}
and the fact that the arithmetic mean curvature lines have
this configuration, as established in \cite{m}.
\end{proof}

\begin{proposition} Consider the  ellipsoid $\mathbb E_{a,b,c}
$ with three axes $a>b>c>0$. On the ellipse $\Sigma\subset \mathbb
E_{a,b,c}$, containing the four  umbilic points,$\;p_i$,
$i=1,\cdots,4,\;$  oriented
counterclockwise, denote by
 $S_1=\int_{-b^2}^{-c^2}\frac{ {1}}{\sqrt{h(v)}}dv$ the distance
between the adjacent umbilic points $p_1 $ and
$p_2 $
and by
$S_2=\int_{-a^2}^{-b^2}\frac{ {1}}{\sqrt{h(u)}  }du$  the   distance
between the adjacent umbilic points
$p_1$ and $p_4$, where $h(x)=(x+a^2)(x+b^2)(x+c^2)$. Define
$\rho=\frac{S_2}{S_1}$.

Then if $\rho\in\Bbb R\setminus \mathbb Q$ (resp. $\rho\in \mathbb
Q$) all the   harmonic mean curvature lines are recurrent ( resp. all,
with the exception of the  harmonic mean
curvature umbilic separatrices, are closed).  \end{proposition}

\begin{proof}
 The ellipsoid $\mathbb E_{a,b,c}$ belongs to  the triple orthogonal
system of surfaces defined by the
one parameter family of quadrics,
$\frac{x^2}{a^2+\lambda}+\frac{y^2}{b^2+
\lambda}+\frac{z^2}{c^2+\lambda}=1$ with $a>b>c>0$, see also
\cite{Sp}  and \cite{St}.

The  following parametrization of $\mathbb E_{a,b,c}$.

$$\alpha(u,v) =
\big(\pm  \sqrt{\frac{M(u,v,a)}{W(a,b,c)}}, \pm
\sqrt{\frac{M(u,v,b)}{W(b,a,c)}}, \pm \sqrt{\frac{M(u,v,c)}{W(c,a,b)}}\big)$$
  \noindent  where,

\noindent $M(u,v,w)= {w^2(u+w^2)(v+w^2)}$ and
$W(a,b,c)=(a^2-b^2)(a^2-c^2)$,
define  the ellipsoidal  coordinates $(u,v)$ on $\mathbb E_{a,b,c}$,
where  $u\in
(-b^2,-c^2)$ and $v\in (-a^2,-b^2)$.

The first fundamental form of $\mathbb E_{a,b,c}$  is given by:

$$I=ds^2=Edu^2+Gdv^2=\frac
14\frac{(u-v)u}{h(u)}du^2+
 \frac 14\frac{(v-u)v}{h(v)}dv^2 $$

The second fundamental form is given by

$$II=
edu^2+gdv^2=\frac{abc(u-v)}{4\sqrt{uv}h(u)
}du^2+\frac{abc(v-u)}{4\sqrt{uv}h(v)}dv^2,$$
\noindent where
$h(x)=(x+a^2)(x+b^2)(x+c^2)$.
  The four umbilic points are $(\pm x_0,0,\pm z_0)= (\pm
a\sqrt{\frac{a^2-b^2}{a^2-c^2}},0,\pm
c\sqrt{\frac{c^2-b^2}{c^2-a^2}}\;)$.

The differential equation of the  harmonic mean curvature lines is
given by:
$$\frac{(du)^2}{h(u)}-\frac{(dv)^2}{h(v)}=0$$

Define $d\sigma_1=  \frac{1}{\sqrt{h(u)}}du$ and
$d\sigma_2=\frac{1}{\sqrt{h(v)}}dv$. By integration, this leads to the
chart $(\sigma_1, \sigma_2)$,
 in which the
differential equation of the  harmonic mean curvature lines is
given by
$$d\sigma_1^2-d\sigma_2^2=0.$$

On the ellipse $\Sigma=\{(x,0,z) | \frac{
x^2}{a^2}+\frac{z^2}{c^2}=1\}$ the distance between the umbilic
points $p_1=(x_0,0,z_0)$ and $p_4=(x_0,0,-z_0)$ is given by
$S_1=\int_{-b^2}^{-c^2}\frac{ {1}}{\sqrt{h(v)}}dv$ and
that between the umbilic points $p_1=(x_0,0,z_0)$ and
$p_2=(-x_0,0,z_0)$
 is given by
$S_2=\int_{-a^2}^{-b^2}\frac{ {1}}{\sqrt{h(u)}  }du$.

It is obvious that  the ellipse $\Sigma$ is the union of four
umbilic points and the   four principal umbilical separatrices for the
principal foliations.  So $\Sigma\backslash\{p_1,p_2,p_3,p_4\}\;$
is a transversal section of both  harmonic mean curvature
foliations. The differential  equation of the  harmonic mean
curvature lines in the principal chart $(u,v)$ is given by
$ e du^2- g dv^2=0$,  which     amounts to $d\sigma_1=\pm
d\sigma_2$. Therefore near the  umbilic point $p_1$ the  harmonic
mean curvature lines with a  harmonic mean curvature umbilic
separatrix contained in the region $\{y>0\}$ define a the transition
map $\sigma_+:\Sigma\to \Sigma$ which is an isometry, reversing
the orientation, with $\sigma_+(p_1)=p_1$. This follows because in
the principal chart $(u,v)$ this
  map is defined by
$\sigma_+:\{u=-b^2\}\to \{v=-b^2\}$ which satisfies the
differential equation $\frac{d\sigma_2}{d\sigma_1}=-1$. By
analytic continuation it results that $\sigma_+$  is an orientation
reversing isometry,
  with two fixed points $\{p_1,\;p_3\}$. The
 harmonic reflection $\sigma_-$, defined in the region $y<0$ have
the two umbilics $\{p_2,\;p_4\}$ as fixed points.

So on the ellipse
parametrized by arclength defined by $\sigma_i$,
  the Poincar\'e return map $\pi_1:\Sigma\to \Sigma$ (
composition of two isometries $\sigma_+$ and $\sigma_-$) is a
rotation with rotation number given by $\frac{S_2}{S_1}$.

Analogously for the other  harmonic mean curvature foliation,
with the Poincar\'e return map given by $\pi_2=\tau_+\circ \tau_-$,
where $\tau_+$ and $\tau_-$ are two isometries having respectively
$\{p_2,p_4\}$ and $\{p_1,p_3\}$ as fixed points. \end{proof}

\section {On Harmonic Mean Curvature Structural Stability}

In this section the results of sections 3, 4 and 5 are put together
 to provide sufficient conditions for harmonic mean curvature
stability, outlined below.

\begin{theorem} \label{th:sta}
The set of immersions ${\mathcal A}_i({\mathbb M^2}), i=1,\; 2$ which
satisfy
conditions $i$), ... , $v$) below
 are i-$C^s$-mean curvature structurally stable and ${\mathcal A}_i,
i=1,\; 2$ is open in ${\mathcal M}^{r,s}({\mathbb M^2)}, \; r\geq  s\geq
6$.
\begin{itemize}

\item[i)]  The parabolic curve is regular : ${\mathcal K}=0$ implies
$d{\mathcal K}   \neq 0$
and the tangential singularities are saddles and nodes.

\item[ii)]  The umbilic points are of type $H_i$, $i=1,\; 2,\;3$.

\item[iii)]  The  harmonic mean curvature cycles of $\mathbb
H_{\alpha,i}$ are hyperbolic.

\item[iv)]  The   harmonic mean curvature foliations $\mathbb
H_{\alpha,i}$  has no separatrix connections. This means that there is no
harmonic mean curvature line joing two umbilic or tangential parabolic
singularities and being separatrices at both ends. See  propositions
\ref{prop:2} and \ref{prop:oc-p}

\item[v)] The limit set of every leaf of $\mathbb H_{\alpha,i}$ is    a
parabolic point, umbilic point or a  harmonic mean curvature cycle.
\end{itemize}

\end{theorem}

\begin{proof} The openness of ${\mathcal A}_i({\mathbb M^2})$ it
follows from the local structure of the  harmonic mean curvature lines near
the umbilic points of types $H_i$, $i=1,2,3$, near the  harmonic mean
curvature cycles and by the absence of umbilic  harmonic mean curvature
separatrix connections and the absence of recurrences.
 The  equivalence can be performed by the method of canonical regions
and their continuation as was done in
\cite{gs1, gs2} for principal lines,  and in \cite{a2}, for asymptotic
lines. \end{proof}

 Notice that   Theorem \ref{th:sta} can be reformulated
so as to give the mean   harmonic stability of the configuration rather
than
that of the  separate  foliations. To this end it is necessary to
consider the folded extended lines, that  is  to consider the line  of one
foliation that arrive  at the parabolic set at a given transversal point
as continuing through the line  of the other foliation leaving the
parabolic set at  this point, in a sort of ``billiard". This gives raise to
the extended folded cycles and separatrices that must be preserved by
the homeomorphism mapping simultaneously the two foliations.

Therefore the third, fourth and fifth hypotheses above should be
modified as follows:
\begin{itemize}
\item[iii')] the extended folded periodic cycles should be hyperbolic,
\item[iv')] the extended folded separatrices should be disjoint,
\item[v')] the limit set of extended lines should be umbilic points,
parabolic singularities and extended folded cycles.
\end{itemize}

 The class of immersions which verify the extended five conditions i),
ii), iii'), iv'), v') of a compact and oriented manifold $\mathbb M^2$
will be denoted by ${\mathcal A}({\mathbb M^2})$.

This procedure has been adopted by the authors in the case of
asymptotic lines by the suspension operation in order to pass from the
foliations to the configuration and properly formulate the stability results.
See \cite{a2}.

\begin{remark} In the space   of convex  immersions ${\mathcal
M}^{r,s}_c({\mathbb S}^2)$ ( ${\mathcal K}_\alpha>0$),  the sets  ${\mathcal
A}({\mathbb S^2})$ and ${\mathcal A}_1({\mathbb S^2})\cap {\mathcal
A}_2({\mathbb S^2}) $ coincide. \end{remark}

 The genericity result involving the five conditions above  is
formulated now.

\begin{theorem} \label{th:gen}
The sets ${\mathcal A}_i, \; i=1,\; 2$ are dense in
${\mathcal M}^{r,2}({\mathbb M^2)},\; r\geq  6$.

In the space  ${\mathcal M}^{r,2}_c({\mathbb S}^2)$ the set ${\mathcal
A}({\mathbb S^2})$ is dense.
\end{theorem}

The main ingredients for the proof of this theorem are the Lifting
and   Stabilization Lemmas, essential for the achievement of
condition five. The conceptual background for this approach goes
back to the works of Peixoto and Pugh.

The elimination of non-trivial recurrences -- the so called  ``Closing
Lemma Problem"-- as a step to achieve condition $v)$is by far the  most
difficult of these details.
See the book of Palis and
Melo, \cite{pm}, for a presentation of these ideas in the case of
vector fields on surfaces.

 The proof of theorem \ref{th:gen} will be postponed to a forthcoming paper \cite {gm}.
It involves   technical details that are closer to
those of the proofs of genericity theorems given by
Gutierrez and Sotomayor, \cite{gsln, gs2}, for principal curvature
lines  and by Garcia and Sotomayor, \cite{m}, for arithmetic  curvature
lines.

\section {Additional Comments and a Related Problem}

The study of families of curves on surfaces  defined by normal
curvature properties and their singularities  has attracted the interest of
generations of mathematicians, among whom can be mentioned Euler,  Monge,
Dupin, Gauss, Cayley, Darboux,  Gullstrand, Caratheodory, Hamburger.
See   \cite{r, St} for references.

On the other hand, the ideas on the {\it``Qualitative Theory of Differential Equations"}
initiated by Poincar\'e and culminating with the study   of the
Structural Stability and Genericity of differential equations on surfaces,
made systematic from $1937$ to $1962$ due  to the seminal work of Andronov
Pontrjagin and Peixoto, were assimilated by Gutierrez, Garcia and
Sotomayor   and, reformulated,  were applied to principal curvature lines
\cite{gs1} as well as to  other  differential equations of classical
geometry: asymptotic lines \cite {a1, a2},  arithmetic and geometric mean
curvature lines \cite {m, g}, and  harmonic mean curvature lines studied
here.

Thus, progress in Differential Equations and Geometry led  to delineate   a
 fruitful field  of interaction
of Geometry and Analysis.

The  work of Monge, on the principal configuration of the  Ellipsoid;
that of  Dupin, on  Triply Orthogonal families of surfaces  and the
study of  Darboux, on   umbilic points on a surface,    are the classical
geometric  paradigms of this field of interaction.

An overview  of the ensemble of  recent  contributions of the authors and others,
cited here,    reveals that there is a common ground. In fact, they share
an  analogy in  purpose, problems  and methods of analysis. It seems,
therefore, appropriate   to inquire here for the common mathematical
features they  enjoy and for the discrepancies they present.

In  principle any expression such as $\mu$= $\mu(k_1,k_2)\in
[k_1,k_2]$, involving  the principal curvatures, could be rightly called a
``mean curvature".

The situations that appear in the works quoted above correspond to  the
Principal Curvatures: $\mu$=$k_1$ or $\mu$=$k_2$,   Arithmetic,
Geometric and Harmonic Mean Curvatures: $\mu$ = $\mathcal H $,
$\mu$ = $\mathcal K^{1/2}$ and  $\mu$ = $\frac{\mathcal
K}{\mathcal H }$.  The asymptotic lines correspond to $\mu=0$. To these five
functions we will refer to as the {\it ``classical" mean curvature
functions.}

\vskip 0.1cm
At this point, a  pertinent problem is proposed  to provoke the
discussion.
\vskip 0.2cm
\noindent {\bf Problem 1.}  Formulate and  prove a  general
theorem from which  the results obtained before for the ``classical"
mean curvature functions  would follow and  that also would include an
interesting  class of ``new" curvature functions $\mu$=
$\mu(k_1,k_2)$ and associated  differential equations.
\vskip 0.2cm

There are  a great number of  {\it means}
that are of interest in Analysis. For instance the
{\it Holder Means}
$$H_r(k_1,k_2)=[(k_1^r+ k_2^r)/2]^{1/r}, $$
which contains,   in the form of a one parameter family,  the classical
means.
In fact,  Arithmetic, corresponds to  $r=1$; Geometric, corresponds to
$r=0$, understood as a limit as ${r\to 0}$;   Harmonic, corresponds to
$r=-1$.  See Hardy et al. \cite {hlp}.

\vskip 0.2cm
 
There are also more subtle, non-algebraic  means,   such as the
AGM-mean,
obtained from the limit of the
Arithmetic and Geometric. This limit was studied classically by
Gauss
and Legendre.  See the book of Borwein and Borwein \cite{bbe}   for the connections of these means with differential
equations and the number $\pi$. The interest of this mean for Geometry
seems to have been overlooked  so far.

A   satisfactory answer to  Problem 1 involves an  
analysis   of the limits of the methods introduced in the recent papers
and their adaptability to deal with configurations associated to
``new" mean curvature functions. Bearing in mind the `` closing
lemma" difficulties mentioned at the end of the previous section,  this
analysis  will be postponed to another work,  \cite{gm}, which contains a proposal for  a partial solution.

\vskip .5cm

\author{\noindent Jorge Sotomayor\\Instituto de Matem\'{a}tica e
Estat\'{\i}stica,\\Universidade de S\~{a}o Paulo, \\Rua do Mat\~{a}o 1010,
Cidade Universit\'{a}ria, \\CEP 05508-090, S\~{a}o Paulo, S.P., Brazil \\
\\ Ronaldo Garcia\\Instituto de Matem\'{a}tica e
Estat\'{\i}stica,\\Universidade Federal de Goi\'as,\\CEP 74001-970, Caixa Postal
131,\\Goi\^ania, GO, Brazil}

\end{document}